\documentclass[11pt]{amsart}
\usepackage{geometry}                
\usepackage{graphicx}
\usepackage{amssymb}
\usepackage{epstopdf}
\usepackage{amsmath}
\usepackage{amsthm}
\usepackage{caption}
\usepackage{subcaption}

\newtheorem{theorem}{Theorem}
\newtheorem{lemma}[theorem]{Lemma}

\theoremstyle{definition}

\newtheorem{xca}[theorem]{Exercise}

\newtheorem{prob}[theorem]{Problem}

\def\vp{\varphi}
\def\al{\alpha}
\def\om{\omega}

\def\fa{\forall}
\def\ex{\exists}
\def\iff{\longleftrightarrow}
\def\then{\longrightarrow}
\def\into{\longrightarrow}

\def\el{elementary}
\def\ele{\el\ end extension}
\def\ar{arithmetic}
\def\ct{countable}
\def\au{automorphism}

\def\cont{$2^{\aleph_0}$}

\def\st{structure}
\def\fo{first-order}

\def\be{\begin{xca}}
\def\ee{\end{xca}}

\def\rs{recursively saturated}
\def\ee{end-extensional}

\def\PA{{\sf PA}}
\def\SPA{{\sf PA}^*}
\def\MS{{\sf MS}}
\def\aca{{\sf ACA}_0}
\def\wkl{{\sf WKL}_0}
\def\act{{\sf ACT}}
\def\rca{{\sf RCA}_0}
\def\srca{{\sf RCA}^*_0}
\def\olike{$\om_1$-like}

\def\zf{{\sf ZF}}

\def\N{{\mathbb{N}}}

\def\X{{\mathfrak{X}}}
\def\Y{{\mathfrak{Y}}}

\DeclareMathOperator{\Th}{Th}

 \DeclareMathOperator{\aut}{Aut}

\DeclareMathOperator{\cf}{cf}  
\DeclareMathOperator{\lt}{Lt}
\DeclareMathOperator{\cod}{Cod}
\DeclareMathOperator{\ssy}{SSy}

\title{63 Years of the MacDowell-Specker Theorem}
\author{ Roman Kossak}

\begin{document}



\maketitle

\section{Introduction} 
In September of 1959,  at the conference on Infinitistic Methods in Warsaw, Ernst Specker presented a joint paper with Robert MacDowell  \cite{MDS}. MacDowell and Specker were studying the additive structure of nonstandard models of \ar. For a brief description of their main result see Hillary Putnam's review  \cite{put}. At the end of the review Putnam writes:

\begin{quotation}
The proof of this theorem requires the result, which is of independent interest, that every model $M$ possesses a proper elementary (``arithmetic") extension $M'$ such that the positive elements of $M'-M$ are greater than the elements of $M$. The proof of this result in turn depends upon an extension of a method due to Skolem. In essence, the method is to form the ultrapower $M^M/{\mathfrak P}$, where $M^M$ is the ring of all (first-order) definable functions from $M$ into $M$ and ${\mathfrak P}$ is a suitable ultrafilter. An ultrafilter with the required properties (the authors speak of a ``finitely additive measure with values 0 and 1'', but this is equivalent to the ultrafilter terminology) is obtained by an ingenious construction. Little previous acquaintance with model theory is required.
\end{quotation}
That result of independent interest is known now as the MacDowell-Specker Theorem. It can be stated simply as follows.
\begin{theorem}[\MS] Every model of Peano Arithmetic has an elementary end extension.\footnote{In this note, elementary extension will always mean proper elementary extension.}
\end{theorem}
MacDowell and Specker proved their theorem for any set of axioms  that contains the Peano axioms (\PA) for addition and multiplication of integers (positive and negative) in the language with constant symbols for 0, 1 and possibly countably many other symbols, with the induction schema for all formulas of the extended language.  
 
In the 60 years that followed, \MS\ has been a constant inspiration for many developments in the model theory of arithmetic. The most striking feature of the theorem is that the same proof covers both the \ct\ and the un\ct\ cases. It is a rarity in model theory. Jim Schmerl has said that every model of \PA\ thinks it is \ct. Thanks to \MS, one can form arbitrarily long chains of \ele s, and one can study unions of such chains that exhibit interesting  properties. Such models  are of independent interest not only in the model theory of \ar. For examples, see  \cite{kauf,  mss, sch85}. 

 \MS\ made the model theory of \PA\ a  special discipline. This becomes evident  when one tries to generalize it to other theories such as fragments of \PA\ or axiomatic systems of set  theory, just to realize that it cannot be done, or at least cannot be done in full generality. Ali Enayat's \cite{ena95} is a relevant survey of corresponding results for models of \zf.

My aim in this note is to discuss briefly the proof of the theorem, and to give an account of its most important  extensions  and generalizations. Then, I will describe the ongoing research and mention some of the main open problems. 

\section{The Theorem}
 
 Using the Omitting Types Theorem, one can show that every {\it countable} model of \PA\ has an elementary end extension. In fact more is true. If a \ct\ structure ${\mathcal M}$ for a countable language is linearly ordered, has no last element  and satisfies the regularity schema consisting of the universal closures of the axioms of the form:
\[ [\fa v \ex x>v\ \ex y<z\ \vp(x,y)]\then [\ex y<z \fa v\ex x>v \ \vp(x,y)],\]
 then ${\mathcal M}$ has an \ele.  This was proved by H.~Jerome Keisler in \cite{kei66, kei71}. A proof for \ct\ models of \PA\  is included  in \cite[Theorem 6.17]{hand}.  For linearly ordered structures  with no last element, the regularity schema is equivalent to the collection principle consisting of the universal closures of the axioms of the form:
 \[ [\fa x< z \ex y_1\dots \ex y_n \vp(x, \bar y)]\then \ex v [\fa x< z \ex y_1<v\dots \ex y_n<v \vp(x, \bar y)].\]
 
For models, $\mathcal M$, $\mathcal N$, \dots, their domains will be denoted $M$, $N$, \dots. For a cardinal number $\kappa$, a linearly ordered model ${\mathcal M}$ is $\kappa$-{\it like} if $|M|=\kappa$ and for each $a$ in $M$, $|\{x: {\mathcal M}\models x<a\}|<\kappa$
 
In  \cite{enamoh},  Enayat and Shahram Mohsenipour  studied the model theory of the regularity schema. They  formulate Keisler's theorem as follows.
 
 \begin{theorem}\label{kei} The following are equivalent for a complete \ct\ \fo\ theory $T$.
 \begin{enumerate}
 \item Some model of $T$ has an \ele.
 \item $T$ proves the regularity schema.
 \item Every \ct\ model of $T$ has an \ele.
 \item Every \ct\ model of $T$ has an \olike\ \ele.
 \item $T$ has a $\kappa$-like model for some regular cardinal $\kappa$.
 \end{enumerate}
 \end{theorem}
 
Because of the limitations on the cardinality of the language in the Omitting Types Theorem, proofs of \MS\ must follow different routes.  MacDowell and Specker modified a construction of Thoralf Skolem. By the compactness theorem,  the standard model of \ar\ $\N=(\om,+,\times,0,1)$ has an \ele.  Skolem proved it directly using what is today called the definable ultrapower construction \cite{sko34}. His proof works for any \ct\ model of \PA. 

Here is a version of Skolem's argument. Let ${\mathcal M}$ be a \ct\ model of \PA, and let $X_n$, $n\in \om$,  be an enumeration of all sets that are definable in ${\mathcal M}$.\footnote{In this note, definability means definability with parameters.} By induction on $n$, we can define a sequence of unbounded definable sets $Y_0\supseteq Y_1\supseteq\cdots$ such that for each $n$, either $Y_n\subseteq X_n$ or  $Y_n\subseteq M\setminus X_n$. Let $p(x)$  be the set of formulas $\vp(x)$ with parameters from $M$, such that for some $n$, $Y_n$ is contained in the set defined by $\vp(x)$ in ${\mathcal M}$. It is easy to check that $p(x)$ is a complete nonprincipal type  of $\Th({\mathcal M},a)_{a\in M}$. 

Using the regularity schema in $\mathcal M$, we can arrange that for each definable $f:M\into M$, if there is an $n$ such that $f(Y_n)$ is bounded in $\mathcal M$, then  $f$ restricted to $Y_m$ is constant for some $m\geq n$.   It is not difficult to check that the model  generated over ${\mathcal M}$ by an element realizing $p(x)$ is an \ele\ of ${\mathcal M}$. 

Using induction in $\mathcal M$, we can  thin down the $Y_n$'s further, so that for every definable $f:M\into M$, there is an $n$ such that  the restriction of $f$ to $Y_n$ is either one-to-one or constant. In this case, the model generated by $M$ and an element realizing  $p(x)$ is a {\it minimal} \ele\ of ${\mathcal M}$, i.e., an \ele\ ${\mathcal N}$ such that ${\mathcal N}$ has no proper elementary submodels properly containing $M$.

 The proof outlined above begins with fixing an enumeration of all formulas of the language of \PA\ with all parameters of a given \ct\ model. The ingenuity of MacDowell and Specker's proof was in their observation that one can avoid having to deal with uncountably many parameters by (almost) forgetting the parameters altogether. Their argument is based on two inductions. One external that applies to an enumeration $\vp_n(x,y)$, $n\in\om$, of all parameter-free formulas of the language of \PA, and one internal, performed in ${\mathcal M}$ for a each $\vp_n(x,y)$ and all instances of $\vp_n(x,b)$, for  $b$ in $M$. For a proof with all  technical details see \cite{kaye1}. 
 
 Here is an abridged proof from \cite{sch73}. Let $\vp_n(x,y)$, $n\in\om$, be a list of all formulas in two variables. Using the fact that induction holds, we can find formulas $\psi_n(y)$ such that for each $n$, the following is a theorem of \PA:
 \[\fa w,z \ex x>z \fa  y<w [\bigwedge_{i<n} (\vp_i(x,y)\iff \psi_i(y))].\]
 Define the type $p(x)$ of  $\Th({\mathcal M},a)_{a\in M}$ by 
 \[ \vp_n(x,b)\in p(x) \textup{ iff } {\mathcal M}\models \psi_n(b).\eqno{(*)}\]
 Again, it  is easy to check that  the model that is generated over ${\mathcal M}$ by an element realizing $p(x)$ is an \ele\ of ${\mathcal M}$. For  details see  \cite[Chapter 3]{ks}.
 
 Once we know that every model of \PA\ has an \ele, there are still many questions that can be asked. An obvious one is: Why is it interesting? I will try to answer it to some extent, but there is also a more technical issue. If we agree, that there is something special about \MS, then we can still ask whether what MacDowell and Specker proved is formulated in the sufficient generality, or, in other words, what is the {\it right} \MS. One way this can be made more precise is as follows.
 
 Let  ${\mathcal N}$ be an \ele, of  ${\mathcal M}$. What can we say about the isomorphism type of the pair $({\mathcal N},M)$? There are two main cases to consider. The first one is when  ${\mathcal N}$ is generated over  ${\mathcal M}$ by a single element, and the second, when it is not the case, and in particular when  ${\mathcal N}$ is the union of a  chain of \ele s.  There is much we could discuss here, but it in this note we will concentrate on the first case. For interesting applications of long chains of MacDowell-Specker extensions see \cite{sch98, sch81a}.

  \section{Conservativity}\label{cons}
The definition $(*)$ of the type $p(x)$ in the previous section turns out to be crucial. A type $p(x)$ of a completion $T$ of \PA\  is {\it definable} if for each parameter-free formula $\vp(x,y)$ there is a formula $\psi(y)$ such that for all constant Skolem terms $t$,
\[ \vp(x,t)\in p(x) \textup{ iff } T\vdash \psi(t).\eqno{(**)}\]
So $(*)$  is an instance of $(**)$ for $T=\Th({\mathcal M},a)_{a\in M}$.

At this point more notation will be useful. Let $p(x)$ be a complete type of $\Th({\mathcal M},a)_{a\in M}$, and let $b$  realize $p(x)$ in ${\mathcal N}$.  Then, let ${\mathcal M}(p)$ be the smallest elementary submodel of ${\mathcal N}$ containing $M\cup \{b\}$. Such a model always exists because \PA\ has definable Skolem functions, and it is unique up to isomorprhism.

If $p(x)$ is a definable type of $\Th({\mathcal M},a)_{a\in M}$, then it follows that for each $X$ that is definable in ${\mathcal M}(p)$, $X\cap M$ is definable in ${\mathcal M}$. Extensions with that property are called {\it conservative}. If ${\mathcal N}$ is a conservative extension of ${\mathcal M}$, then for every $c\in N\setminus M$, $\{x\in M: {\mathcal N}\models x<c\}$ is definable in ${\mathcal M}$. This implies that every conservative extension is an end extension. MacDowell and Specker did not mention conservativity, but a careful reading of their proof reveals that every model of \PA\ has a conservative \ele.  Robert Phillips defined conservativity in \cite{phi74}, and  gave an explicit proof that every model of \PA\ has a conservative \ele. 

The {\it standard system}, $\ssy({\mathcal M})$,  of a nonstandard model ${\mathcal M}$ is the set of the sets of the form $\om\cap X$, for all $X$ that are definable in ${\mathcal M}$. By compactness, for every set of natural numbers $X$, the standard model $\N$ has an elementary extension ${\mathcal M}$ such that $X$ is in $\ssy({\mathcal M})$. Hence, $\N$ has  nonconservative extensions. 

Phillips proved in \cite{phi74a} that the standard model has nonconservative minimal \ele s. This can be generalized to all \ct\ models of \PA\ as follows.  By a result of Stephen Simpson \cite{sim74}, every \ct\ model of ${\mathcal M}$ of \PA\ has undefinable subsets $X$ such that $({\mathcal M},X)\models \SPA$, i.e., $({\mathcal M},X)$ satisfies the induction schema in the language with an additional unary predicate symbol interpreted as $X$. By \MS, $({\mathcal M},X)$ has an \ele\ $({\mathcal N},Y)$. By induction on $a$, one can show that for each $a$ in $N$, the set $Y\cap \{x: N\models x<a\}$ is coded by an element of $N$, hence it is definable in ${\mathcal N}$. For each $a$ in $N\setminus M$, $X=M\cap (Y\cap \{x: N\models x<a\})$; hence ${\mathcal N}$ is not a conservative extension of ${\mathcal M}$.  If ${\mathcal N}$ is an elementary extension of ${\mathcal M}$, then a subset $X$ of $M$ is {\it coded} in ${\mathcal N}$, if for some $Y$ that is definable in ${\mathcal N}$, $X=M\cap Y$.  It is shown in \cite{kospar} that every subset of a \ct\ model of $\PA$ that can be coded in an \ele, can be coded in a minimal \ele.

The set of subsets of a model ${\mathcal M}$ that are coded in an \ele\ ${\mathcal N}$ is denoted by $\cod({\mathcal N}/{\mathcal M})$.
In Section \ref{coded}, we will discuss in detail the question: Which sets subsets of a model ${\mathcal M}$ are of the form $\cod({\mathcal N}/{\mathcal M})$ for some \ele\ of ${\mathcal M}$?

A subset $X$ of a model ${\mathcal M}$ of \PA\ is a {\it class} if for each $a$ in $M$, $\{a: {\mathcal M}\models x<a\}\cap X$ is definable in ${\mathcal M}$. Each definable set is a class, and so is each set $X$ such the $({\mathcal M},X)\models\SPA$. Also, if ${\mathcal N}$ is an \ele\ of ${\mathcal M}$, then all sets in $\cod({\mathcal N}/{\mathcal M})$ are classes of ${\mathcal M}$. Schmerl took the term ``rather classless" from the title of \cite{kauf} and  defined in \cite{sch81a} a model to be {\it rather classless} if all of its classes are definable.

Conservativity in \MS\ is not an incidental property. By \MS, every model of \PA\ has arbitrarily long chains of  conservative \ele s. Schmerl showed  that if $\cf(\kappa)>\aleph_0$, and ${\mathcal N}$ is $\kappa$-like model of \PA\, that is the union of a $\kappa$-chain  of conservative \ele s, then ${\mathcal N}$ is rather classless  \cite[Theorem 1.6]{sch73}. All \ele\ of rather classless models are conservative. This suggests that every proof of \MS\ must, explicitly or implicitly, yield the existence of conservative extensions. 

\section{Elementarity}
In \cite{MDS}, MacDowell and Specker used \MS\ to show that for every model ${\mathcal M}$ of \PA\ and every cardinal $\kappa\geq |M|$, there  is a model of $\PA$, ${\mathcal N}$ such that  $|N|=\kappa$ and $\ssy({\mathcal M})=\ssy({\mathcal N})$. This is a simple consequence of the fact that, by \MS, every model of \PA\ has arbitrarily long chains of \ele s. 

To prove that every model of \PA\ has an {\it end extension}, one can apply the Arithmetized Completeness Theorem (\act). \act\ implies that every model of \PA\ has a conservative end extension. The construction can be iterated along $\om$, but then it breaks down. The union of a chain of end extensions of models of \PA\ may not be a model of \PA.  Kenneth McAloon constructed  $\om$-chains of end extensions of models of \PA\ such that in their unions $\om$ is definable by a $\Delta^0_1$ formula \cite{mca78}.

Attempts to generalize \MS\ to fragments of \PA\ have led to interesting problems. If a model of $I\Delta^0_0$ has an elementary end extension, then it is a model of \PA. This was shown by Laurence Kirby and Jeff Paris in \cite{parkir}, where they also proved the following theorem.

\begin{theorem}\label{clote} For every $n\geq 2$, every \ct\ model of $B\Sigma^0_n$ has a $\Sigma^0_n$-elementary end extension.
\end{theorem}
Peter Clote tried to show that Theorem \ref{clote} holds for uncountable models as well, but only succeeded to do it with the stronger assumption that ${\mathcal M}\models I\Sigma^0_n$,  \cite{clo98, clo86}. The problem whether the result holds for $B\Sigma^0_n$ is open.  For a discussion of the end extension problem for fragments of \PA, see  \cite{dim}, where, among other results, Costas Dimitracopoulos gives an alternative proof of Clote's theorem using \act.

If $\kappa$ is a regular cardinal and ${\mathcal M}$ is a $\kappa$-like model of $I\Delta^0_0$, then, by the results of Kirby and Paris, ${\mathcal M}$ is a model of \PA.  For singular $\kappa$ the situation is more complex. Richard Kaye proved in \cite{kay97} that  for each $n\geq 1$ and  each singular $\kappa$, there are $\kappa$-like models of $B\Sigma^0_n+\exp+\lnot I\Sigma^0_n$, and he asked 
whether  $B\Sigma^0_1$ + exp+ $\{I\Sigma^0_n\then B\Sigma^0_{n+1}: n\in\om\}$ is an axiomatization of the theory of all $\kappa$-like models for singular strong limit cardinals $\kappa$.\footnote{$I\Sigma_n$ and $B\Sigma_{n+1}$ are finitely axiomatized for all $n\geq 1$.} This was answered in the negative by Ian Robert Haken and Theodore Slaman, who constructed a $\kappa$-like model of $I\Sigma^0_1+\lnot B\Sigma^0_2$ for a singular strong limit $\kappa$, \cite[Chapter 3]{haksla}.

\section{Canonicity}
The study of elementary extensions of models of \PA\ was raised to another level by Gaifman's seminal paper \cite{gai76}. 

Let $T$ be a completion of \PA. A type $p(x)$ is  {\it unbounded} if it contains all formulas $(t<x)$, where $t$ is a constant Skolem term; $p(x)$ is {\it end-extensional} if it is unbounded, and for every model ${\mathcal M}$ of $T$, if ${\mathcal N}$ is generated over ${\mathcal M}$ by an element $b$ realizing $p(x)$, and $a<b$ for all $a$ in $M$, then ${\mathcal N}$ is an \ele\ of ${\mathcal M}$. If $p(x)$ is unbounded and each ${\mathcal N}$ as above is a minimal \ele\ of ${\mathcal M}$, then $p(x)$ is  {\it minimal}.

By the definition, each minimal type is \ee, and Gaifman shows that each \ee\ type is definable, and that none of these inclusions reverses. Moreover,  the set of \ee\ types is dense in the space of complete unbounded 1-types of $T$.  

There are several equivalent definitions of minimal types (see \cite[Theorem 3.2.10]{ks}). In particular, each minimal type $p(x)$ is {\it strongly indiscernible}, i.e., for every model ${\mathcal M}$ of \PA, for every set of elements $I$  realizing $p(x)$ in ${\mathcal M}$, and every $a$ in $I$, for all $c_0, c_1,\dots, c_k\leq a$, $\{x\in I: a<x\}$ is a set of indiscernibles in the structure $({\mathcal M},c_0,c_1,\dots, c_k)$. 

If $p(x)$ is an \ee\ type of $T$ and ${\mathcal M}$ and ${\mathcal N}$ are isomorphic models of $T$, then ${\mathcal M}(p)$ and ${\mathcal N}(p)$ are isomorphic as well. Moreover, in the standard topological space of countable models, the operation ${\mathcal M}\mapsto {\mathcal M}(p)$ is Borel. 

We can also define a canonical operation ${\mathcal M}\mapsto {\mathcal M}(I)$, where $(I,<)$ is a linearly ordered set, and $M(I)$ is generated over $M$ by a set of indiscernibles $\{a_i: i\in I\}$, realizing  a fixed  minimal type and such that ${\mathcal M}(I)\models a_i<a_j$ iff $i<j$. Gaifman proved that for all ${\mathcal M}$ and $(I,<)$, the group of \au s of ${\mathcal M}(I)$ that fix ${\mathcal M}$ pointwise is isomorphic to $\aut(I,<)$ \cite[Theorem 4.11]{gai76}. This is improved by the following theorem of Jim Schmerl.

\begin{theorem}[\cite{sch02}] For every linearly ordered \st\ $\mathfrak A$, every  model of \PA\ has an \ele\ ${\mathcal N}$ such that $\aut({\mathcal N})\cong \aut({\mathfrak A})$.
\end{theorem}

The operation ${\mathcal M}\mapsto {\mathcal M}(I)$, reduces the isomorphism relation for \ct\ linear orders  to the isomorphism relation of \ct\ models of $T$. This shows that the latter relation is Borel complete. Along with the MacDowell-Specker-Gaifman techniques, there are other ways of constructing \ele s with special properties. Canonicity of those constructions has not been systematically studied, but there is  one interesting example which we will discuss next.

We say that ${\mathcal N}$ is a {\it superminimal} extension of ${\mathcal M}$ if for every $a\in N\setminus M$, the only elementary submodel of ${\mathcal N}$ containing $a$ is ${\mathcal N}$. Every \ct\ model of $\PA$ has a conservative superminimal \ele. The proof given in \cite[Theorem 2.1.12]{ks} begins with enumerating $M$, and it is not a good prognosis for  canonicity of the operation \[{\mathcal M}\mapsto\textup{ a superminimal \ele\ of } {\mathcal M}.\]
It is shown in \cite{ck} that in fact there is no  operation that would in a Borel way map isomorphic models of \PA\ to their isomorphic superminimal \ele s.

If ${\mathcal N}$ is the  union of a chain ${\mathcal M}_\al$, $\al<\om_1$, of models such that for each $\al$, ${\mathcal M}_{\al+1}$ is a superminimal \ele\ of ${\mathcal M}_\al$, then ${\mathcal N}$ is {\it J\'onsson}, i.e., ${\mathcal N}$ has no proper elementary submodel of cardinality $\aleph_1$. By a theorem of Andrzej Ehrenfeucht \cite{ehr73}, ${\mathcal N}$ realizes $\aleph_1$ complete types.

\begin{prob} Is there an un\ct\ J\'onsson model of \PA\  that realizes only $\aleph_0$  complete types?
\end{prob}

Gaifman defines a dependence relation on types as follows. Let $p(x)$ and $q(x)$ be complete types of a completion of \PA. Then, $q(x)$ {\it depends} on $p(x)$ if there is a Skolem term $t(x)$ such that for all $\vp(x)$, $\vp(x)\in q(x)$ iff $\vp(t(x))\in p(x)$. Gaifman showed that for every completion $T$, there are \cont\  independent minimal types of $T$ \cite[Theorem 4.13]{gai76}. It follows that every \ct\ model of \PA\ has \cont\  conservative minimal \ele s.

Schmerl has shown that  there is a family of minimal types $\{p_\sigma(x): \sigma \in 2^{\om}\}$, such that for any $\sigma$, $\tau$ in $2^\om$, $p_\sigma(x)$ and $p_\tau(x)$ are dependent if and only if there is a $k$ such that for all $n\geq k$, $\sigma(n)=\tau(n)$ \cite[Lemma 3.6]{ck}. This result was used in \cite{ck} to show that for any completion $T$, the isomorphism relation for finitely generated models of $T$ is not Borel reducible to the  identity relation on $2^\om$. 

In \cite{sh78},  Saharon Shelah looked for possible  extensions of Gaifman's result for theories other than \PA. He wrote
\begin{quotation}
Now Gaifman \cite{gai65, gai67, gai76}, following MacDowell and Specker \cite{MDS}, proved that any theory $T=\Th(\om,+,\times,<,\dots)$ has definable end extensions, minimal ones, rigid ones, etc. He uses the fact that {\it definitions by induction} are allowable. We show that for some of the results {\it proof by induction} only is sufficient, so $T=\Th(\om,<,\dots)$  is sufficient (for almost minimal end-extension types). This seems maximum we can get, but we do not have counterexamples.
\end{quotation}

Shelah called a model ${\mathcal N}$  an {\it almost minimal} \ele\ extension of ${\mathcal M}$ if all elementary submodels of ${\mathcal N}$ that contain $M$ are cofinal in ${\mathcal N}$, and he defined {\it almost  minimal \ee} types. Shelah proved that every  theory $T$  in a \ct\ language that proves the induction schema and the sentence $\fa x \ex y\ (x<y)$, has almost minimal \ee\ types.

In \cite{ena07}, Enayat intoduced a general framework  for constructions based on iterations  of  what he calls {\it Skolem-Gaifman} ultrapowers, and  other types of restricted ultrapowers over models of \PA, and he  obtains several  interesting results on automorphism of \ct\ \rs\ models of \PA.

\section{Uncountable languages}
\MS\ and Gaifman's results  hold for all theories that include Peano's axioms in \ct\ languages extending   the language of \PA. Therefore, many later results on models of \PA\ extend without modifications to such languages. To indicate that, they are formulated for $\SPA$, without specifying the language.

Gaifman  asked if the main results of \cite{gai76} can be proved for extensions of \PA\ in uncountable languages. That quickly turned out not to be the case. George Mills used a version of arithmetic forcing to show that  each countable model ${\mathcal M}$ of \PA\ can be expanded by  a set of functions $f_\al: M\into M$, $\al<\om_1$, so that the model $({\mathcal M},f_\al)_{\al<\om_1}$ satisfies Peano's axioms, and it has no \ele\ \cite{mil78}. Nevertheless, some positive results can be proved as well. 

In \cite{sch73}, Schmerl defined generic sets of subsets of a model of \PA, and he proved the following generalization of \MS\ (proofs are included in \cite[Chapter 6]{ks}): 
\begin{theorem}[ \cite{sch81,sch73}] Let ${\mathcal M}$ be a model of \PA.
\begin{enumerate}
\item For every generic $\mathfrak X$,  $({\mathcal M},A)_{A\in {\mathfrak X}}$ satisfies the induction schema.

\item For every generic $\mathfrak X$,  $({\mathcal M},A)_{A\in {\mathfrak X}}$ has a conservative \ele.

\item For any infinite cardinal $\kappa$, there are models of \PA\  of cardinality $\kappa$ with generic sets of subsets of cardinality $\kappa^+$
\end{enumerate}
\end{theorem}

In \cite{bla74}, Andreas Blass extended some of Gaifman's results to models of \PA\ in the full language of \ar, i.e., the language with function and relation symbols for all functions and relations of the standard model. Blass characterized \ee\ and minimal types of this theory in terms of ultrafilters on $\om$. 

Every elementary extension of the standard model is an end extension, but even in this case \MS\ gives additional information: the standard model has a conservative elementary end extension. 

 Every extension of $(\N,A)_{A\subseteq \om}$ is  conservative, but there are expansions of the standard model $\N$ without elementary conservative extensions.  This was shown by Enayat in \cite{ali08}. The problem, listed as Question $7'$ in \cite[Chapter 12]{ks}, whether there is an expansion $\mathfrak N$  of the standard model such that some nonstandard model of $\Th({\mathfrak N})$ has no (conservative) \ele,  was settled in the positive direction by Shelah in \cite{she11}.

\section{The Lattice Problem}
For a model ${\mathcal M}$ of \PA, the set of elementary submodels of ${\mathcal M}$ forms a lattice $\lt({\mathcal M})$, with the  join of two models being the smallest elementary submodel of ${\mathcal M}$ containing both of them, and the meet being their intersection. The general lattice problem is to characterize those lattices that can be represented as $\lt({\mathcal M})$, for some  ${\mathcal M}$. While much is known about the problem, it is still open whether every finite lattice can be represented this way. For a given ${\mathcal M}$, one can also ask about lattices of interstructures $\lt({\mathcal N}/{\mathcal M})$, where ${\mathcal N}$ is an elementary extension of ${\mathcal M}$, and $\lt({\mathcal N}/{\mathcal M})$ is the set of elementary submodels of ${\mathcal N}$ that contain $M$. In particular, we can ask about lattices that can be realized as $\lt({\mathcal N}/{\mathcal M})$, where ${\mathcal N}$ is an \ele\ of ${\mathcal M}$. Here Gaifman's technique comes to play.

Let $T$ be a completion of $\PA$ and let $p(x)$  be an \ee\ type of $T$. By a theorem of Gaifman \cite{gai76}, for every model ${\mathcal M}$ of $T$, $p(x)$ can be (uniquely) extended to a complete definable type $p'(x)$ of $\Th({\mathcal M},a)_{a\in M}$, such that $a<x$ is in $p'(x)$ for all $a$ in $M$.  Since $p'(x)$ is uniquely determined by $p(x)$, we will write  ${\mathcal M}(p)$ instead of ${\mathcal M}(p')$. 

After Mills \cite{mil78}, we will say a type of $T$ {\it produces} a lattice $L$, if for every model ${\mathcal M}$ of $T$, $\lt({\mathcal M}(p)/{\mathcal M})$ is isomorphic to $L$. Thus,  every  minimal $p(x)$ produces the two-element Boolean algebra $\{{\bf 0,1}\}$. Gaifman proved several results on lattices representable as substructure and  interstructure lattices,  and  conjectured that for every finite distributive lattice $L$ with exactly one atom there is an \ee\ type that produces $L$ \cite{gai76}. The conjecture was confirmed by Schmerl \cite{sch78},
and extended further by Mills \cite{mil79}.  Mills extended Gaifman's technique of \ee\ 1-types to types with arbitrary (finite or infinite) sets of variables, and he characterized completely distributive lattices (of any cardinality) that can be produced by \ee\ types. 

\begin{figure}
\centering
\begin{subfigure}{.5\textwidth}
  \centering
  \includegraphics[width=.4\linewidth]{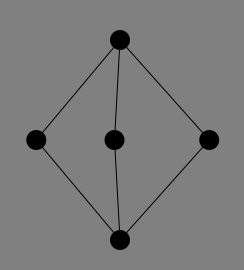}
  \caption{}
  \label{fig1}
\end{subfigure}%
\begin{subfigure}{.5\textwidth}
  \centering
  \includegraphics[width=.4\linewidth]{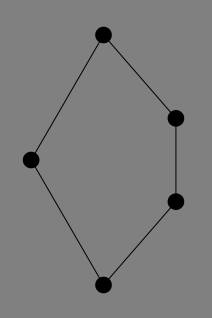}
  \caption{}
  \label{fig2}
\end{subfigure}
\caption{${\bf M}_3$ and ${\bf N}_5$.}
\label{fig}
\end{figure}

The case of the distributive lattices is closed, but there is still much mystery about the nondistributive ones. While, as shown in \cite{sch86},  the lattice ${\bf M}_3$ (Figure \ref{fig},A) can be represented as $\lt({\mathcal N}/{\mathcal M})$, where ${\mathcal N}$ is a cofinal extension of ${\mathcal M}$, ${\bf M}_3$ has no interstructure lattice representation in which ${\mathcal N}$ is an \ele\ of ${\mathcal M}$ \cite{gai76, par77}. Alex Wilkie proved in \cite{wil77} that  every \ct\ model ${\mathcal M}$ has an \ele\ such that $\lt({\mathcal N}/{\mathcal M})$ is isomorphic to ${\bf N}_5$ (Figure \ref{fig},B). 

It is not easy to say why the cases of ${\bf M}_3$, and ${\bf N}_5$ are so different, and what makes them both different from the distributive lattices case. All constructions of \ele s require some Ramsey style combinatorics, nontrivial results in lattice representation theory, and proofs involve some rather intricate details. 
Much  work on the lattice problem after the 1980's was done by Schmerl who has developed special techniques for constructing models with prescribed properties. An account of those developments can be found in \cite[Chapter 4]{ks}. One of the results which is directly related to the theme of this note that is not covered there is the following theorem.
\begin{theorem}[\cite{sch04}] Let $T$ be a completion of \PA, and let $L$ be a finite lattice.
\begin{enumerate}
\item If some countable model ${\mathcal M}$ of $T$ has an almost minimal \ele\ ${\mathcal N}$ such that $\lt({\mathcal N}/{\mathcal M})\cong L$, then every countable model ${\mathcal M}$ of $T$ has an almost minimal \ele\ ${\mathcal N}$ such that $\lt({\mathcal N}/{\mathcal M})\cong L$.
\item If some countable model ${\mathcal M}$ of $\Th(\N)$ has an almost minimal \ele\ ${\mathcal N}$ such that $\lt({\mathcal N}/{\mathcal M})\cong L$, then every  model ${\mathcal M}$ of $\Th(\N)$ has an almost minimal \ele\ ${\mathcal N}$ such that $\lt({\mathcal N}/{\mathcal M})\cong L$.
\end{enumerate}
\end{theorem}

 The next problem is listed as  Question 4 in \cite[Chapter 12]{ks}.
\begin{prob} What finite lattices $L$ are such that every ${\mathcal M}$ has an \ele\ such that $\lt({\mathcal N}/{\mathcal M})\cong L$? What finite lattices $L$ are such that every {\it \ct} ${\mathcal M}$ has an \ele\ such that $\lt({\mathcal N}/{\mathcal M})\cong L$? What finite lattices $L$ are such that {\it some} \ct\ ${\mathcal M}$ has an \ele\ such that $\lt({\mathcal N}/{\mathcal M})\cong L$?
\end{prob}

\begin{prob} Let $L$ be a finite lattice, and suppose that each model ${\mathcal M}$ of \PA, has an \ele\ ${\mathcal N}$ such that $\lt({\mathcal N}/{\mathcal M})\cong L$. Is there an \ee\ type $p(x)$ such that for every ${\mathcal M}$, $\lt({\mathcal M}(p)/{\mathcal M})\cong L$?
\end{prob}

\section{Coded sets}\label{coded}
Recall that if ${\mathcal N}$ is an \ele\ of ${\mathcal M}$, then  $\cod({\mathcal N}/{\mathcal M})$ denotes  the set of subsets of $M$ that are coded in  ${\mathcal N}$. With this notation, we can rephrase \MS\ as follows: 

\begin{theorem} Every model ${\mathcal M}$ of $\SPA$ has an \ele\ ${\mathcal N}$ such that $\cod({\mathcal N}/{\mathcal M})$ is the collection of all definable subsets of $M$. 
\end{theorem} 
Let $(\N,{\mathfrak X})$ be a  model of $\aca$, and let $({\mathcal M},A^*)_{A\in {\mathfrak X}}$ be a conservative \ele\ of $(\N,A)_{A\in {\mathfrak X}}$. Then $\X=\ssy({\mathcal M})$. Hence, for each \ct\ model $(\N,{\mathfrak X})$ of $\aca$ there is an \ele\ ${\mathcal M}$ of $\N$ such that $\X=\ssy({\mathcal M})$. This is special case of Scott's theorem that characterizes the \ct\ $\om$-models of $\wkl$ as standard systems of \ct\ models of \PA, \cite{sco62}. The same argument  gives the following result.

\begin{theorem}\label{cor} If ${\mathcal M}$ and $\mathfrak X$ are \ct\ and  $({\mathcal M},{\mathfrak X})$ is a model of $\aca$, then ${\mathcal M}$ has an \ele\ such that $\cod({\mathcal N}/{\mathcal M})={\mathfrak X}$.
\end{theorem}

Theorem \ref{cor} was generalized by Schmerl in the next theorem which I am tempted to call the ultimate \MS. It is not difficult to show that if $\mathcal N$ is an \ele\ of $\mathcal M$, then $\cod({\mathcal N}/{\mathcal M})$ contains all definable subsets of $\mathcal M$, and that $({\mathcal M}, \cod({\mathcal N}/ {\mathcal M}))$ is a model of $\wkl^*$, which is the usual $\wkl$ of second-order \ar\ in which $I\Sigma^0_1$ is replaced with $I\Delta^0_0+\exp$.  A set $\X$ of subsets of $M$ is {\it countably generated} if there is a \ct\ $\X_0$ such that every set in $\X$ is $\Delta^0_1$-definable in $({\mathcal M},\X_0)$.  

\begin{theorem} [\cite{sch14}]\label{sch14} Let $\X$ be a set of subsets of a model ${\mathcal M}$ of \PA. Then the following are equivalent: 
\begin{enumerate} 
\item $\X$ is countably generated,  it contains all definable subsets of ${\mathcal M}$, and $({\mathcal M},\X)$ is a model of $\wkl^*$. 
\item  ${\mathcal M}$ has a finitely generated elementary end extension ${\mathcal N}$ of ${\mathcal M}$ such that $\X=\cod({\mathcal N}/{\mathcal M})$. 
\item  ${\mathcal M}$ has a countably generated elementary end extension ${\mathcal N}$ of ${\mathcal M}$ such that $\X=\cod({\mathcal N}/{\mathcal M})$. 
\end{enumerate}
\end{theorem}

 The question whether this Theorem \ref{sch14} can be strengthened further by requiring that ${\mathcal N}$ is a minimal extension was also settled by Schmerl, with an intriguing twist:
 
\begin{theorem}[\cite{sch17}]\label{sch17} Let $\X$ be a set of subsets of a model ${\mathcal M}$ of \PA. Then the following are equivalent: 
\begin{enumerate} 
\item  $\X$ is countably generated,  it contains all definable subsets of ${\mathcal M}$,  $({\mathcal M},\X)$ is a model of $\wkl^*$, and every set that is $\Pi^0_1$-definable in $({\mathcal M},\X)$ is the union  of countably many $\Sigma^0_1$-definable sets.
\item ${\mathcal M}$ has a minimal \ele\ ${\mathcal N}$ such that $\X=\cod({\mathcal N}/{\mathcal M})$.
 \end{enumerate}
 \end{theorem}
 
 Proofs of both theorems are very finely tuned versions of the original proof of MacDowell and Specker. 
 
 It is not clear what the rather exotic condition in Theorem \ref{sch17} (1) really means, but, by the next theorem, we know that it is necessary. Notice that if  $({\mathcal M},\X)\models \aca$, then $\X$ trivially satisfies the condition.

 \begin{theorem}[\cite{sch17}] For every model ${\mathcal M}_0$ of \PA, there are ${\mathcal M}$ and $\X$ such that $\Th({\mathcal M}_0)=\Th({\mathcal M})$, $({\mathcal M},\X)\models\wkl$, $\X$ is countably generated and contains all definable subsets of ${\mathcal M}$, and there is a set that is $\Pi^0_1$-definable in $({\mathcal M},\X)$ that is not the union of countably many $\Sigma^0_1$-definable sets.
 \end{theorem}

 Athar Abdul-Quader has observed (unpublished) that Theorem \ref{sch17} can be extended by adding the third condition: For each finite distributive lattice $D$, ${\mathcal M}$ has an \ele\ ${\mathcal N}$ such that $\X=\cod({\mathcal N}/{\mathcal M})$ and $\lt({\mathcal N}/{\mathcal M})\cong D$. Schmerl gave a short proof of this result using arguments from \cite[Chapter 4]{ks}.

 Another interesting development was triggered by an attempt to generalize Wilkie's theorem about ${\bf N}_5$ to the un\ct\ case.
 It turns out that this can't be done in general. If $\lt({\mathcal N}/{\mathcal M})$ is isomorphic to ${\bf N}_5$, then ${\mathcal N}$ is not a conservative extension. An interesting proof, due to Schmerl, is included in \cite[Section 4.6]{ks}. It follows that if $\mathcal M$ is rather classless and $\mathcal N$ is a \ele, then $\lt({\mathcal N}/{\mathcal M})$ is not isomorphic to  ${\bf N}_5$.  This implies that, ${\bf N}_5$  cannot be produced by an \ee\ type.   
 
 Matt Kaufmann has observed that \MS\ can be  proved using the Arithmetized Completeness Theorem. An outline of his proof is 
 given in \cite{sch92}. In \cite{enawon}, Enayat and Tin Lok Wong use a refined version of this argument to give an alternative proof of a slightly weaker version of Theorem \ref{sch14}.  It is not clear if their method yields finitely generated extensions.

 \begin{prob} Let  ${\mathcal M}$ be a model of \PA.
 \begin{enumerate}
 \item If  $\lt({\mathcal N}/{\mathcal M})\cong {\bf N}_5$, what can we say about $\cod({\mathcal N}/{\mathcal M})$?
 \item If ${\mathcal M}$ is \ct, is it true that for every undefinable subset  $A$ of ${\mathcal M}$ there is an \ele\ ${\mathcal N}$ such that $\lt({\mathcal N}/{\mathcal M})\cong {\bf N}_5$ and $A$ is not in $\cod({\mathcal N}/{\mathcal M})$?
 \end{enumerate}
 \end{prob}

In the recent paper \cite{simwon}, Simpson and Tin Lok Wong combine several results from the literature to give the following characterizations of $\wkl$  and  $\aca$. The extension $(N,\Y)$ of $(M,\X)$ is {\it conservative} if for all $Y\in\Y$, $Y\cap M\in \X$.
 
\begin{theorem} Let  $(M,\X)$ be a \ct\ model of $\srca$. Then the following are equivalent:
\begin{enumerate}
\item $(M,\X)\models\wkl.$
\item $(M,\X)\models I\Sigma^0_1$ and it has a proper conservative extension that satisfies $\srca$.
\end{enumerate}
\end{theorem}

\begin{theorem} Let  $(M,\X)$ be a  \ct\ model of $\srca$. Then the following are equivalent:
\begin{enumerate}
\item $(M,\X)\models\aca.$
\item $(M,\X)$ has a proper $\Sigma^1_1$-elementary conservative extension that satisfies $\aca$.
\end{enumerate}
\end{theorem}
In \cite{sch86a}, Schmerl proved that every \ct\ model of $\aca$ has a $\Sigma^1_1$-elementary {\it exclusive}\footnote{See \cite{sch86a} or \cite{simwon} for a definition.} extension which satisfies $\aca$. Simpson and Tin Lok Wong give an alternative proof of Schmerl's result, and they characterized \ct\ models of ${\sf ATR}_0$ and $\Pi^1_1$-${\sf CA}_0$ in terms of exclusive extensions.

Kirby proved that a \ct\ model of $\rca$ admits a nonprincipal definable type if and only if it satisfies $\aca$ \cite{kir84}. Simpson and Tin Lok Wong pose the following problem.

 \begin{prob} Can one characterize \ct\ models of ${\sf ATR}_0$ and $\Pi^1_1$-${\sf CA}_0$ in terms of some variants of definable types, without requiring the defining scheme to be in any specific form?
 \end{prob}
 \section{Recursive saturation}
 
A systematic study of the model theory of \ct\ \rs\ models of \PA\ was initiated by Craig Smory\'nski in \cite{smo}. Smory\'nski proved, among many other results,  that every \ct\ \rs\ model ${\mathcal M}$ of \PA\ has an  \ele\ that is isomorphic to ${\mathcal M}$. Prior to that, John Schlipf \cite{schl78} proved that resplendent model of \PA\ has a resplendent \ele. So  there is a version of \MS\ for resplendent models,  but  it does not extend to all  \rs\  models. There are \rs\ models of \PA\ that have no \rs\ \ele. The prime example is the \olike\ rather classless \rs\ model constructed by Kaufmann in \cite{kauf}.  It is not difficult to show that if ${\mathcal N}$ is a \rs\ \ele\ of ${\mathcal M}$, then the extension is not conservative. This implies that Kaufmann's model has no \rs\ \ele s. Kaufman's proof uses Jensen's diamond principle which can be eliminated  by  Shelah's Absoluteness Theorem \cite{she78}. Still, the following problem, posed by Wilfrid Hodges in \cite{hod85}, is open.

\begin{prob}
Prove the existence of rather classless \rs\ models of \PA\ in cardinality $\aleph_1$ without assuming diamond at any stage of the proof.
\end{prob}

Here is the key lemma in Kaufman's proof:

\begin{lemma}\label{kauf} Let ${\mathcal M}$ be a \ct\ \rs\ model of \PA. If $A$ is an undefinable subset of $M$, then ${\mathcal M}$ has a \rs\ \ele\ ${\mathcal N}$ such that 
$A\notin\cod({\mathcal N}/{\mathcal M})$.
\end{lemma}

In \cite{sch81a}, Schmerl studied rather classless \rs\ model of \PA\ in cardinalities higher than $\aleph_1$ by means of special types of MacDowell-Specker extensions. Here is one of his main results.
\begin{theorem} Assume $V=L$. Let $\kappa$ be an infinite cardinal and let $T$ be a completion of \PA. Then the following are equivalent:
\begin{enumerate}
\item there is  a $\kappa$-like, \rs\ model of $T$;
\item $\cf(\kappa)>\aleph_0$ and $\kappa$ is not weakly compact.
\end{enumerate}\end{theorem} 

Lemma \ref{kauf} says that each undefinable subset of $M$ can  be omitted from $\cod({\mathcal N}/{\mathcal M})$ in some \rs\ \ele\ of ${\mathcal M}$.  The following problem asks if in a similar manner we can also omit the complete theory of any undefinable subset of $M$. It is listed as Problem 17 in \cite[Chapter 12]{ks}, but is badly marred by typos. Here is the correct version.

\begin{prob}\label{kaufp} Let $A$ be an undefinable subset of a \ct\ \rs\ model ${\mathcal M}$ of \PA. Does ${\mathcal M}$ have a  \rs\ \ele\ ${\mathcal N}$ such that for every $B\in\cod({\mathcal N}/{\mathcal M})$, $\Th({\mathcal M},A)\not=\Th({\mathcal M},B)$?
\end{prob}

If $\Th({\mathcal M},A)$ is not in $\ssy({\mathcal M})$, then the answer to the above question is positive. By chronic resplendency,   ${\mathcal M}$ has a \rs\ \ele\ ${\mathcal N}$ such that $({\mathcal N},{\mathcal M})$ is \rs. In this case, for each $B\in\cod({\mathcal N}/{\mathcal M})$, $\Th({\mathcal M},B)\in\ssy({\mathcal M})$. 

There is an open problem concerning \ele\ of \rs\ models that is similar in flavor to Problem \ref{kaufp}. It was first posed in \cite{kos4}, and despite some effort, it remains open. 

\begin{prob}
Let ${\mathcal M}$ be \ct\ and \rs\ and let ${\mathcal K}$ be a proper cofinal submodel of ${\mathcal M}$. Does ${\mathcal M}$ have a \rs\ \ele\ ${\mathcal N}$ such that for every \rs\  ${\mathcal N}'$ such that  ${\mathcal K}\preceq {\mathcal N}'\prec {\mathcal N}$, if $N'\cap M=K$, then $N'=K$.
\end{prob}

\bibliographystyle{plain}
\bibliography{MDS}

\end{document}